\documentclass[a4paper,11pt]{amsart}

\usepackage{amssymb}
\usepackage{amscd}
\usepackage[all]{xy}

\def\frak{\mathfrak}
\def\Bbb{\mathbb}
\def\Cal{\mathcal}

\newtheorem*{prop*}{Proposition}

\newtheorem*{thm*}{Theorem}

\newtheorem*{lem*}{Lemma}

\newtheorem*{cor*}{Corollary}

\newcommand{\id}{\operatorname{id}}
\newcommand{\im}{\operatorname{im}}
\newcommand{\tr}{\operatorname{tr}}
\newcommand{\End}{\operatorname{End}}

\newcommand{\x}{\times}
\renewcommand{\o}{\circ}

\let\ccdot\cdot
\def\cdot{\hbox to 2.5pt{\hss$\ccdot$\hss}}

\newcommand{\ka}{\kappa}
\newcommand{\la}{\lambda}
\newcommand{\om}{\omega}
\renewcommand{\phi}{\varphi}

\newcommand{\ps}{\psi}

\newcommand{\ze}{\zeta}

\newcommand{\La}{\Lambda}
\newcommand{\Ga}{\Gamma}
\newcommand{\Ph}{\Phi}
\newcommand{\Om}{\Omega}

\newcommand{\Rho}{{\mbox{\sf P}}}

\renewcommand{\k}{{\mathbf{k}}}

\newcommand{\ce}{{\Cal E}}

\def\sideremark#1{\ifvmode\leavevmode\fi\vadjust{\vbox to0pt{\vss
 \hbox to 0pt{\hskip\hsize\hskip1em
 \vbox{\hsize3cm\tiny\raggedright\pretolerance10000
 \noindent #1\hfill}\hss}\vbox to8pt{\vfil}\vss}}}%
\begin{document}
\title[A holonomy characterisation of Fefferman spaces]{A holonomy characterisation\\ of Fefferman spaces}
\author{Andreas \v Cap and A.\ Rod Gover} 
\date{November 30, 2006}

\address{A.C.: Fakult\"at f\"ur Mathematik, Universit\"at Wien,
  Nordbergstra\ss e 15, A--1090 Wien, Austria and International Erwin
  Schr\"odinger Institute for Mathematical Physics, Boltzmanngasse 9,
  A--1090 Wien, Austria\newline\indent A.R.G.: Department of
  Mathematics, The University of Auckland, Auckland, New Zealand}
\email{Andreas.Cap@esi.ac.at, gover@math.auckland.ac.nz}
\begin{abstract}
  We prove that Fefferman spaces, associated to non--de\-ge\-ne\-rate CR
  structures of hypersurface type, are characterised, up to local
  conformal isometry, by the existence of a parallel orthogonal
  complex structure on the standard tractor bundle. This condition can
  be equivalently expressed in terms of conformal holonomy. Extracting
  from this picture the essential consequences at the level of tensor
  bundles yields an improved, conformally invariant analogue of
  Sparling's characterisation of Fefferman spaces.
 \end{abstract}

\maketitle

\section{Introduction}\label{1}
Fefferman spaces provide a geometric relationship between CR geometry
and conformal geometry; given a CR manifold $M$, of hypersurface type,
one obtains a canonical conformal structure on the total space $\tilde
M$ of a certain circle bundle over $M$. The first version of this
construction in \cite{Fefferman:annals} applied to the boundaries of
strictly pseudoconvex domains, and used Fefferman's ambient metric
construction.  Soon after that, a version for abstract non--degenerate
CR structures of hypersurface type was given in \cite{BDS}. This used
the canonical Cartan connection for CR structures of Chern and Moser
\cite{Chern-Moser}.

Recently there has been renewed interest in the Cartan connection
approach to CR structures and, more generally, parabolic
geometries. Several powerful tools, like tractor calculus and BGG
sequences for these geometries have been developed. These tools
provide a new and effective approach the construction of Fefferman
spaces and to treating the complicated relationship between the
natural objects on the Fefferman space with those on the underlying
CR manifold. This approach has been taken up in the article
\cite{fefferman}, where several new results on Fefferman spaces their
conformal geometry are obtained. Central in these developments is the
result that on the (conformal) standard tractor bundle of a Fefferman
space $\tilde M$, one obtains a parallel orthogonal complex
structure. This complex structure makes the standard tractor bundle
into a Hermitian vector bundle and the standard tractor connection
into a Hermitian connection. Key to the power of this approach is that
these data are rather simply related to the CR standard tractor bundle
of $M$.

In this article we complete the picture by showing that Fefferman
spaces are in fact characterised by this orthogonal parallel complex
structure, up to local conformal isometry. This has a very nice
interpretation in terms of conformal holonomy, an area of significant
recent interest, see e.g. \cite{Armstrong} and references therein. By
definition, conformal holonomy is the holonomy group of the standard
tractor connection. In the dimensions and signatures which are
relevant for us, this is a subgroup of $SO(2p'+2,2q'+2)$. Existence of
a parallel orthogonal complex structure is clearly equivalent to
the conformal holonomy being contained in the subgroup
$U(p'+1,q'+1)$. Hence we obtain a characterisation of Fefferman spaces
among general conformal structures analogous to the characterisation
of K\"ahler manifolds among pseudo-Riemannian manifolds.  By
construction, the conformal holonomy of a Fefferman space is in fact
automatically contained in the smaller subgroup $SU(p'+1,q'+1)$. As a
corollary, we therefore obtain that if the conformal holonomy is
contained in $U(p'+1,q'+1)$ then the holonomy Lie algebra already has
to be contained in $\frak{su}(p'+1,q'+1)$. In particular,
$\frak{u}(p'+1,q'+1)$ cannot be a conformal holonomy Lie algebra
although it is a Berger algebra. In fact this follows easily and
directly from a simple calculation using the tractor connection, see
Proposition \ref{2.1}. This was proved in a more involved way in
\cite{Leitner}.

 There is a characterisation of Fefferman spaces available in the
literature, which is usually referred to as Sparling's
characterisation. While Sparling's original work remained unpublished,
the characterisation has been proved in a different way by C.R.~Graham
in \cite{Graham}. The central ingredient for this characterisation is
an isotropic Killing field with certain additional properties. In our
picture, this arises as follows: A parallel orthogonal complex
structure $\Bbb J$ on the standard tractor bundle can be interpreted
as a parallel section of the adjoint tractor bundle. Any section of
the adjoint tractor bundle has an underlying vector field, and for a
parallel section this is automatically a conformal Killing field which
inserts trivially into the Cartan curvature. Expanding the condition
$\Bbb J\o\Bbb J=-\id$ yields the remaining properties required in
Sparling's characterisation (and other identities that are
differential consequences of these).

In the last section of the article we that show that the tensorial
consequences of $\Bbb J\o\Bbb J=-\id$ lead to a characterisation of
Fefferman spaces that is in the same spirit as Sparling's, but is
conformally invariant. Using this, with the extra choice of a certain
conformal scale, recovers Sparling's characterisation precisely. This
provides additional insight into the structure of conformal manifolds
which admit isotropic conformal Killing fields.

We assume that the reader is familiar with basic conformal geometry
and conformal tractor calculus and we will follow the conventions of
\cite{conf-Einst}. Other sources for background on conformal tractor
calculus are \cite{BEG} and \cite{confamb}, for generalisations to
parabolic geometries see \cite{TAMS}. All facts about the construction
of Fefferman spaces we need can be found in \cite{fefferman}.

\subsection*{Acknowledgement} First author supported by Project
P15747-N05 of the Fonds zur F\"orderung der wissenschaftlichen
Forschung (FWF). The second author would like to thank the Royal
Society of New Zealand for support via Marsden Grant no.\ 02-UOA-108,
and the New Zealand Institute of Mathematics and its Applications for
support via a Maclaurin Fellowship.

\section{A characterisation of Fefferman spaces}\label{2}

\subsection{Parallel adjoint tractors and conformal Killing 
  fields}\label{2.1} 

The main technical input for the characterisation is provided by
general results in conformal geometry. Hence we start in the general
setting of a smooth manifold $M$ of dimension $n\geq 3$ endowed with a
conformal structure $[g]$ of signature $(p,q)$. Then the
\textit{standard tractor bundle} (see \cite{BEG,confamb}) is a vector
bundle $\Cal T\to M$ of rank $n+2$, endowed with a bundle metric $h$
of signature $(p+1,q+1)$, a line subbundle $\Cal T^1\subset\Cal T$
with isotropic fibres, and a linear connection $\nabla^{\Cal T}$ which
is compatible with $h$ and satisfies a certain non--degeneracy
condition with respect to the subbundle $\Cal T^1$. These data are
canonically associated with the conformal structure and are equivalent
to its canonical Cartan connection.

The \textit{adjoint tractor bundle} $\Cal A\to M$ is then defined as
the bundle $\frak{so(\Cal T)}$ of skew symmetric endomorphisms of
$\Cal T$. The adjoint tractor connection $\nabla^{\Cal A}$ is the
linear connection on $\Cal A$ induced by $\nabla^{\Cal T}$. There is a
canonical projection $\Pi:\Cal A\to TM$ from the adjoint tractor
bundle to the tangent bundle, so in particular, any section of $\Cal
A$ has an underlying vector field.

We will sometimes work with abstract indices, denoting tractor indices
by $A,B,\dots$ and tensor indices by $a,b,\dots$. Tractor indices are
raised and lowered using the tractor metric $h=h_{AB}$ and its
inverse, while tensor indices are raised and lowered using the
\textit{conformal} metric ${\bf g}_{ab}\in\Ga(\Cal E_{ab}[2])$ and its
inverse. Thus raising (lowering) a tensor index decreases (increases)
the conformal weight by $2$. If we use connections in abstract index
computations, then we will always assume that a metric from the
conformal class has been chosen. The symbol $\nabla_a$ then denotes a
coupled Levi--Civita--tractor connection, i.e.~it acts by the standard
tractor connection on all tractor indices and by the Levi--Civita
connection on all tensor indices.

The curvature $\Om$ of the standard tractor connection $\nabla^{\Cal
  T}$ can be naturally interpreted as an element of $\Om^2(M,\Cal A)$.
Its abstract index representation therefore has the form
$\Om_{ab}{}^A{}_B$. Lowering the tractor index $A$, the result is
skew symmetric in both pairs of indices. The form $\Om$ also describes
the curvature of the canonical Cartan connection. It can be easily
expressed explicitly (see the proof below) in terms of the
\textit{Weyl--curvature} $C=C_{ab}{}^c{}_d$ and the
\textit{Cotton--York tensor} $A_{abc}$.

\begin{prop*}
    Let $s=s^A{}_B$ be a parallel section of the adjoint tractor
    bundle $\Cal A$ and let $\k=\k^a$ denote the underlying vector
    field $\Pi(s)$. Then we have:

\noindent
(1) $\k$ is a conformal Killing field which in addition satisfies
    $\k^a\Om_{ab}{}^A{}_B=0$; 

\noindent
(2) $\Om_{ab}{}^A{}_Bs^B{}_A=0$.
  \end{prop*}
  \begin{proof}
    (1) is proved in \cite[Proposition 2.2]{conf-Einst} and in a much
    more general setting in \cite[Corollary 3.5]{deformations}. 

The explicit formula for the tractor curvature from \cite[formula
(6)]{conf-Einst} reads as 
$$
\Om_{abAB}=Z_A{}^cZ_B{}^dC_{abcd}-2X_{[A}Z_{B]}{}^cA_{cab},
$$
with the convention that the last two indices of the Cotton tensor
are the two--form indices. Hence from (1) we get $\k^aC_{abcd}=0$ and
$\k^a A_{cab}=0$. From Proposition 2.2 and Lemma 2.1 of
\cite{conf-Einst} we get   $s_{AB}=\tfrac{1}{n}D_AK_B$, where
$K_B=Z_B{}^a\k_a-\tfrac{1}{n}X_B\nabla_a\k^a$ for any choice of metric
in the conformal class. Since $K_B$ has conformal weight one, we
obtain
$$
s_{AB}=Y_AK_B+Z_A{}^a\nabla_a K_B-
\tfrac{1}{n}X_A(\nabla^a\nabla_a+\Rho^a{}_a)K_B.
$$
Contracting both tractor indices with $\Om_{abAB}$, the last
summand does not contribute, since \mbox{$\Om_{abAB}X^A=0$}, and a
short computation gives 
$$
\Om_{abAB}s^{AB}=-2\k^cA_{cab}+C_{abcd}\nabla^c\k^d .
$$
 
For the first term, we observe that the total alternation of the
Cotton--York tensor vanishes. (This is a consequence of symmetry of
the Schouten tensor $\Rho_{ab}$.) Hence
$$
\k^cA_{cab}=\k^c(-A_{bca}-A_{abc}),
$$
and we have observed above that the right hand side vanishes.
Finally, the symmetries of the Weyl tensor together with
$\k^bC_{abcd}=0$ give
$$
C_{abcd}\nabla^c\k^d=-\k^d\nabla^cC_{abcd}=-(n-3)\k^dA_{dab}=0,
$$
where we have used the Bianchi identity in the last step. This
completes the proof of (2).
\end{proof}

\subsection{Complex structures on standard tractors}\label{2.2}
Let us now assume that $\dim(M)$ is even and $[g]$ is a conformal
class of signature $(2p'+1,2q'+1)$, so the tractor metric has
signature $(2p'+2,2q'+2)$. We will further assume that we have given
an orthogonal complex structure $\Bbb J$ on the standard tractor
bundle $\Cal T$, which is possible for such a signature. Observe that
since $\Bbb J^2=\Bbb J\o\Bbb J=-\id$, orthogonality of $\Bbb J$ is
equivalent to skew--symmetry. Hence $\Bbb J$ can be considered as a
section of the adjoint tractor bundle $\Cal A$. Note also that, using
$\Bbb J$, we can extend the tractor metric $h$ to a Hermitian bundle
metric $\Cal H$.

By definition, $\Bbb J$ is parallel for $\nabla^\Cal A$ if and only if
covariant derivatives by $\nabla^\Cal T$ are complex linear. In
this case $\nabla^\Cal T$ is a Hermitian
connection for $\Cal H$ and hence its holonomy is contained in
$U(p'+1,q'+1)\subset SO(2p'+2,2q'+2)$. Conversely, this condition on
the holonomy is evidently equivalent to existence of a parallel
orthogonal complex structure $\Bbb J$. 

Using $\Bbb J$, we can view the standard tractor bundle $\Cal T$ as a
complex vector bundle. In particular, we can form the complex line
bundle $\Cal V:=\La^{p'+q'+2}_{\Bbb C}\Cal T$, i.e.~the highest complex
exterior power of $\Cal T$. The standard tractor connection
$\nabla^\Cal T$ induces a linear connection $\nabla^\Cal V$ on $\Cal V$, 
and the
curvature of this linear connection is given by the trace of the
tractor curvature, viewed as a complex linear map. More explicitly,
this is the complex valued two--form
$\Om_{ab}{}^A{}_A-i\Om_{ab}{}^A{}_B\Bbb J^B{}_A$. The first summand
vanishes by skew symmetry while the second vanishes by part (2) of
Proposition \ref{2.1}. Therefore, locally the bundle $\Cal V$ admits
smooth sections which are parallel for the linear connection induced
by $\nabla^{\Cal T}$.

We want to interpret the existence of $\Bbb J$ via a reduction of
structure group of the canonical Cartan bundle $\Cal G$ associated to
the conformal structure. The Cartan bundle can be easily obtained from
the standard tractor bundle $\Cal T$ by a frame bundle construction,
see \cite[2.2]{confamb}. To prepare the grounds for using complex
structures later, we use the underlying real vector space of $\Bbb
V:=\Bbb C^{p'+q'+2}$ as the modelling vector space for $\Cal T$. Fix a
Hermitian inner product $\langle\ ,\ \rangle$ of signature
$(p'+1,q'+1)$ on $\Bbb V$ as well as a real isotropic line $L_\Bbb
R\subset\Bbb V$, and a non-zero element $\nu_0\in\La^{p'+q'+2}_\Bbb
C\Bbb V$. Let $\langle\ ,\ \rangle_{\Bbb R}$ be the real part of
$\langle\ ,\ \rangle$. Let $\Cal T^1\subset\Cal T$ be the
distinguished real line subbundle in $\Cal T$. Then the fibre $\Cal
G_x$ of the canonical Cartan bundle can be identified with the set of
all orthogonal linear isomorphisms $u:\Bbb V_\Bbb R\to\Cal T_x$ such
that $u(L_{\Bbb R})=\Cal T^1_x$. Fixing an element $u\in\Cal G_x$, we
obtain a bijection from the stabiliser $P\subset SO(\Bbb V)=:G$ of the
line $L_\Bbb R$ onto $\Cal G_x$.  Together with the fact that $\Cal T$
can be locally trivialised as a filtered metric vector bundle, this
implies that the disjoint union $\Cal G=\sqcup_{x\in M}\Cal G_x$ can
be made into a smooth principal $P$--bundle over $M$. By construction,
$\Cal T$ is the associated bundle $\Cal G\x_P\Bbb V$.

Let $\frak g=\frak{so}(\Bbb V)$ be the Lie algebra of $G$ and let
$\frak p\subset\frak g$ be the subalgebra corresponding to $P$. The
canonical conformal Cartan connection $\om\in\Om^1(\Cal G,\frak g)$
can be recovered from the standard tractor connection as follows.
Since $\Cal T=\Cal G\x_P\Bbb V$, the space $\Ga(\Cal T)$ of smooth
sections can be naturally identified with the space $C^\infty(\Cal
G,\Bbb V)^P$ of smooth functions $f:\Cal G\to\Bbb V$ such that
$f(u\cdot g)=g^{-1}(f(u))$ for all $u\in\Cal G$, and all $g\in
P\subset SO(\Bbb V)$. Any vector field $\xi\in\frak X(M)$ can be
lifted (locally) to a vector field $\overline{\xi}$ on $\Cal G$. For a
section $s\in\Ga(\Cal T)$ corresponding to $f:\Cal G\to\Bbb V$, the
function inducing $\nabla_\xi^\Cal T s$ is then given by
$\overline{\xi}\cdot f+\om(\overline{\xi})\o f$, for any such lift
$\xi$.  It is shown in \cite[2.7]{TAMS} that this uniquely defines a
Cartan connection $\om$ on $\Cal G$ and normality of $\om$ is
equivalent to normality of $\nabla^\Cal T$.

Now suppose that $\Bbb J$ is an orthogonal complex structure on $\Cal
T$. From above we then know that $\Cal V=\La^{p'+q'+2}_{\Bbb C}\Cal T$
admits local non--vanishing parallel sections. Restricting to an open
subset if necessary, we assume that $\nu$ is a global nonzero parallel
section of $\Cal V$. Then we define $\Cal G^{SU}_x\subset\Cal G_x$ to
be the set of those maps $u$, which are complex linear with respect to
$\Bbb J_x$, and which also have the property that the induced
isomorphism on the $(p'+q'+2)$nd complex exterior power maps $\nu_0$
to $\nu(x)$. The subset of complex linear maps in $SO(\Bbb V)$ is
exactly the unitary group $U(\Bbb V)$, such a map lies in $SU(\Bbb
V)$, if in addition the induced map on $\La^{p'+q'+2}_{\Bbb C}\Bbb V$
preserves $\nu_0$.  Thus, $\Cal G^{SU}:=\sqcup\Cal G^{SU}_x\subset\Cal
G$ is a principal fibre bundle over $M$ with structure group
$P^{SU}:=SU(\Bbb V)\cap P$.  This structure group is exactly the
stabiliser of the real line $L_{\Bbb R}$ in the group $SU(\Bbb V)$.

Note that by elementary linear algebra $SU(\Bbb V)$ acts transitively
on the space of real null lines, so the inclusion $SU(\Bbb
V)\hookrightarrow G=SO(\Bbb V)$ induces a diffeomorphism $SU(\Bbb
V)/P^{SU}\cong G/P$. Looking at derivatives at the base points, we
see that the inclusion $\frak{su}(\Bbb V)\hookrightarrow\frak g=\frak
{so}(\Bbb V)$ induces a linear isomorphism $\frak{su}(\Bbb V)/\frak
p^{SU}\to \frak g/\frak p$.

\begin{prop*}
  The canonical conformal Cartan connection $\om$ restricts to a
  Cartan connection $\om^{SU}\in\Om^1(\Cal G^{SU},\frak{su}(\Bbb V))$.
\end{prop*}
\begin{proof}
  The relationship between $s\in\Ga(\Cal T)$ and the corresponding
  equivariant function $f\in C^\infty(\Cal G,\Bbb V)^P$ is
  characterised by the fact that for each $u\in\Cal G_x$ we get
  $s(x)=u(f(u))\in\Cal T_x$. Elements of $\Cal G^{SU}$ are by
  definition complex linear as isomorphisms $\Bbb V\to\Cal T_x$.
  Consequently, if $\tilde f:\Cal G\to\Bbb V$ is the equivariant
  function corresponding to $\Bbb Js$, then $\tilde f(u)=if(u)$ for
  all $u\in\Cal G^{SU}$. Now assume that $\xi\in\frak X(M)$ is a
  vector field, and consider a local lift $\overline{\xi}$ on $\Cal
  G^{SU}\subset\Cal G$, which is tangent to $\Cal G^{SU}$. Then the
  integral curves of $\tilde\xi$ are contained in $\Cal G^{SU}$, so
  along these integral curves we have $\tilde f=if$. Thus we have
  $\overline{\xi}\cdot\tilde f=i(\overline{\xi}\cdot f)$, and thus
  $\nabla^\Cal T_\xi \Bbb Js$ corresponds to $i(\overline{\xi}\cdot
  f)+\om(\overline{\xi})\o if$.  On the other hand, $\nabla^\Cal T_\xi
  \Bbb Js=\Bbb J\nabla^\Cal T_\xi s$, which, along $\Cal G^{SU}$
  corresponds to the function $i(\overline{\xi}\cdot
  f+\om(\overline{\xi})\o f)$. But this exactly shows that
  $\om(\overline{\xi})$ is complex linear, provided that
  $\overline{\xi}$ is tangent to $\Cal G^{SU}$.
  
  According to the observations above, we can view $\Cal T$ as the
  associated bundle $\Cal G^{SU}\x_{P^{SU}}\Bbb V$. Thus the exterior
  power $\Cal V$ is the associated bundle $\Cal
  G^{SU}\x_{P^{SU}}\La^{p'+q'+2}\Bbb V$.  Now elements of $SU(\Bbb V)$
  by definition act trivially on this exterior power, so $\Cal V$ is
  actually the trivial bundle $M\x\Bbb C$. By construction, the
  section $\nu$ of $\Cal V$ corresponds to the constant function
  $\nu_0$. The relation between $\nabla^\Cal V$ and $\om$ is similar
  to the situation of $\nabla^\Cal T$, as discussed above. Using this,
  we can conclude that for $\overline{\xi}$ tangent to $\Cal G^{SU}$,
  we have $\om(\overline{\xi})\in\frak{su}(\Bbb V)$.
  
  Hence we can restrict $\om$ to $\om^{SU}\in\Om^1(\Cal
  G^{SU},\frak{su}(\Bbb V))$. On each tangent space, the map
  $\om^{SU}$ is evidently injective and hence bijective for
  dimensional reasons.  Denoting by $r^g:\Cal G\to\Cal G$ the principal
  right action of $g\in P$, equivariancy of $\om$ reads as
  $(r^g)^*\om=\operatorname{Ad}(g^{-1})\o\om$. For $g\in P^{SU}\subset
  P$ we get $r^g(\Cal G^{SU})\subset\Cal G^{SU}$ and equivariancy of
  $\om$ immediately implies that
  $(r^g)^*\om^{SU}=\operatorname{Ad}(g^{-1})\o\om^{SU}$ for such $g$.
  Finally, take an element $A$ in the Lie algebra $\frak p^{SU}$ of
  $P^{SU}$.  Viewing $A$ as an element of $\frak p$, we get the
  fundamental vector field $\ze_A\in\frak X(\Cal G)$ generated by $A$.
  By definition, in points of $\Cal G^{SU}\subset\Cal G$, the field
  $\ze_A$ is tangent to $\Cal G^{SU}$. Hence we can restrict $\ze_A$
  to a vector field on $\Cal G^{SU}$ and by definition, this
  restriction coincides with the fundamental vector field generated by
  $A$ on this smaller principal bundle. Thus $\om(\ze_A)=A$ implies
  that $\om^{SU}$ reproduces the generators of fundamental vector
  fields. This completes the verification that $\om^{SU}$ defines a
  Cartan connection on $\Cal G^{SU}$.
\end{proof}

\subsection{Passing to a local leaf space}\label{2.4}
By part (1) of Proposition \ref{2.1}, the vector field $\k^a$,
underlying the parallel adjoint tractor $\Bbb J$, is a conformal
Killing field. The argument of \cite[Theorem 3.1]{fefferman} shows
that $\k^a$ is nowhere vanishing and by \cite[4.4]{fefferman} it is
actually a Killing field for appropriate metrics in the conformal
class. Since $\k$ is nowhere vanishing, it defines a rank one
foliation of $M$, and we may consider a local leaf space for this
foliation. This means that for each $x\in M$ we find an open
neighbourhood $W$ of $x$ and a smooth surjective submersion $\ps:W\to
N$ onto some manifold $N$, such that $\ker(T_x\ps)=\Bbb R\k(x)$ for
each $x\in W$.

A particular outcome from Proposition \ref{2.2} is that the tangent
bundle $TM$ can be naturally viewed as the associated bundle $\Cal
G^{SU}\x_{P^{SU}}\frak{su}(\Bbb V)/\frak p^{SU}$, where $\frak p^{SU}$
is the Lie algebra of $P^{SU}$. In this picture, there is a helpful
interpretation of the distribution spanned by $\k$. It is represented
by the real span of the class of any element $A\in\frak{su}(\Bbb V)$
such that $A-\Bbb J\in\frak p\subset\frak{so}(\Bbb V)$. But the fact
that $A-t\Bbb J\in\frak p$ for some $t\in\Bbb R$ is equivalent to the
fact that $A$ stabilises the complex line $L\subset\Bbb V$, which is
generated by $L_{\Bbb R}$.

Let us denote by $\frak q\subset\frak{su}(\Bbb V)$ and $Q\subset
SU(\Bbb V)$ the stabilisers of $L$. Then $Q$ is a parabolic subgroup
of $SU(\Bbb V)$ with Lie algebra $\frak q$. By construction, $Q$ acts
on the line $L$ and the stabiliser of the real line $L_\Bbb R\subset
L$ is the subgroup $P^{SU}\subset Q$.  Since the action of $Q$ on $L$
is evidently transitive, we see that $Q/P^{SU}\cong\Bbb RP^1$, so in
particular this quotient is connected.

\begin{thm*}
  Let $(M,[g])$ be a conformal manifold endowed with a parallel
  orthogonal complex structure $\Bbb J$ on the standard tractor bundle
  $\Cal T$ and a non vanishing parallel section $\nu$ of $\Cal V$. Let
  $(p:\Cal G^{SU}\to M,\om^{SU})$ be the Cartan geometry of type
  $(SU(\Bbb V),P^{SU})$ obtained in Proposition \ref{2.2}.
  
  Then, for sufficiently small local leaf spaces $\ps:W\to N$, we can
  find a Cartan connection $\underline{\om}\in\Om^1(N\x
  Q,\frak{su}(\Bbb V))$ on the trivial principal bundle $N\x Q$ and a
  $P^{SU}$--equivariant diffeomorphism $\Ph$ from a
  $P^{SU}$--invariant open subset of $N\x Q$ onto a
  $P^{SU}$--invariant subset of $\Cal G^{SU}$ which pulls back
  $\om^{SU}$ to $\underline{\om}$.
\end{thm*}
\begin{proof}
  We follow the proofs of Proposition 2.6 and Theorem 2.7 of
  \cite{twistors}, and refer to the notation introduced just above the
  Theorem. By part (1) of Proposition \ref{2.1}, $\k$ hooks trivially
  into the Cartan curvature, so the mapping
  $A\mapsto(\om^{SU})^{-1}(A)\in\frak X(\Cal G^{SU})$ is a Lie algebra
  homomorphism. This can be viewed as an action of $\frak q$ on $\Cal
  G^{SU}$, and such an action integrates to a local action of $Q$ on
  $\Cal G^{SU}$. Choosing the leaf space $N$ so small that there is a
  local section $s$ of $\Cal G^{SU}\to W\to N$ one finds an open
  neighbourhood $V$ of $P^{SU}$ in $Q$ which is invariant under right
  multiplication by elements of $P^{SU}$. Using $s$, one constructs a
  $P^{SU}$--equivariant diffeomorphism $\Ph$ from $N\x V$ onto a
  $P^{SU}$--invariant subset of $\Cal G^{SU}$, such that $\ps\o
  p\o\Ph=\text{pr}_1:N\x Q\to N$ and for $A\in\frak q$, with
  corresponding fundamental vector field $\ze_A\in\frak X(N\x Q)$, we
  have $T\Ph\o\ze_A=(\om^{SU})^{-1}(A)\o\Ph$.
  
  The fullback $\Ph^*\om^{SU}\in\Om^1(N\x V,\frak{su}(\Bbb V))$
  restricts to a linear isomorphism on each tangent space. The
  restriction of this pullback to $N\x\{e\}$ can be extended
  equivariantly to a Cartan connection $\underline{\om}\in\Om^1(N\x
  Q,\frak{su}(\Bbb V))$. Using that $\k$ inserts trivially into the
  Cartan curvature, one shows that $\underline{\om}$ coincides with
  $\Ph^*\om^{SU}$ locally around $N\x\{e\}$ and hence on $N\x V$ by
  equivariancy.
\end{proof}

\begin{cor*}
  A local leaf space $N$ as in the Theorem inherits an almost CR
  structure of hypersurface type. Explicitly, the complex subbundle
  $H\subset TN$ of corank one is the image of the orthocomplement
  $\k^\perp\subset TM$.
\end{cor*}
\begin{proof}
  The Cartan geometry $(N\x Q,\underline{\om})$ gives rise to an
  identification of $TN$ with the associated bundle $(N\x
  Q)\x_Q(\frak{su}(\Bbb V)/\frak q)\cong N\x (\frak{su}(\Bbb V)/\frak
  q)$. By definition, $\frak q$ is the stabiliser of the complex line
  $L\subset\Bbb V$. Now we define $\frak h\subset \frak{su}(\Bbb
  V)/\frak q$ as the set of those $A+\frak q$ for which $A(L)\subset
  L^\perp$, where $L^\perp$ is the complex orthocomplement of $L$.
  This is a well defined subspace in $\frak{su}(\Bbb V)/\frak q$.
  Fixing a nonzero element $v\in L_{\Bbb R}$ the map $A+\frak q\mapsto
  A(v)+L$ induces a linear isomorphism $\frak{su}(\Bbb V)/\frak q\to
  v^{\perp_{\Bbb R}}/L$, where we use the real orthocomplement of $v$.
  This restricts to an isomorphism $\frak h\cong L^\perp/L$, so $\frak
  h\subset \frak{su}(\Bbb V)/\frak q$ has real codimension one.
  Moreover, $L^\perp/L$ is a complex vector space, and the induced
  complex structure on $\frak h$ is independent of the choice of $v$.
  Since both the subspace $\frak h$ and its complex structure are
  invariant under the natural action of $Q$, they induce a complex
  subbundle $H\subset TN$ of real codimension one.
  
  To describe $H$ explicitly, consider a point $u\in\Cal
  G^{SU}\subset\Cal G$ and a tangent vector $\xi\in T_u\Cal G^{SU}$.
  By definition $T\ps\cdot Tp\cdot\xi$ lies in the CR subbundle if and
  only if $A:=\om(\xi)\in\frak {su}(\Bbb V)$ maps $L$ to $L^\perp$.
  Let $v\in L_{\Bbb R}$ be a nonzero element. Since $A$ is orthogonal,
  we have $\langle A(v),v\rangle_{\Bbb R}=0$, and using that $A$ is
  complex linear we see that $A(L)\subset L^\perp$ if and only if
  $\langle A(v),iv\rangle_{\Bbb R}=0$. But, by construction, on the
  subbundle $\Cal G^{SU}$, the element $i\id+\frak p\in\frak {so}(\Bbb
  V)/\frak p$ corresponds to the vector field $\k$. Since the
  conformal structure is induced by $\langle\ ,\ \rangle_{\Bbb R}$,
  the explicit description of $H$ follows.
\end{proof}

\subsection{The curvature of the reduced Cartan connection}\label{2.5}
The next step is to show that the almost CR structure from Corollary
\ref{2.4} is non--degenerate of signature $(p',q')$ and integrable. We
will prove more than that, namely that the Cartan geometry $(N\x
Q,\underline{\om})$ is torsion free and normal, which also implies
that it is the canonical Cartan geometry associated to this CR
structure. To study these issues, we have to understand the curvature
of $\underline{\om}$.

For a sufficiently small local leaf space $\ps:W\to N$, Theorem
\ref{2.4} provides us with a Cartan geometry of type $(SU(\Bbb V),Q)$
on $N$. Regular Cartan geometries of that type induce an underlying
partially integrable almost CR structure, and conversely, a partially
integrable almost CR structure gives rise to a unique regular normal
Cartan connection, see \cite[2.3]{fefferman} and
\cite[4.15]{Cap-Schichl}.

Recall that the curvature of the Cartan connection $\om\in\Om^1(\Cal
G,\frak g)$ is the two form $K\in\Om^2(\Cal G,\frak g)$ defined by
$K(\xi,\eta)=d\om(\xi,\eta)+[\om(\xi),\om(\eta)]$ for
$\xi,\eta\in\frak X(\Cal G)$. The defining properties of a Cartan
connection immediately imply that this form is horizontal and and
$P$--equivariant. Using the trivialisation of $T\Cal G$ provided by
$\om$, one obtains the curvature function $\ka:\Cal G\to L\big(\La^2(\frak
g/\frak p),\frak g\big)$, which is characterised by
$$
\ka(u)(X+\frak p,Y+\frak p):=K(\om^{-1}(X)(u),\om^{-1}(Y)(u)).
$$
This is well defined since $K$ is horizontal and equivariancy of
$K$ implies that $\ka$ is equivariant for the natural $P$--action on
$L\big(\La^2(\frak g/\frak p),\frak g\big)$. On the other hand, composing the
inverse of the isomorphism $\frak{su}(\Bbb V)/\frak p^{SU}\to\frak
g/\frak p$ induced by inclusion, with a natural projection, we obtain a
natural surjection $\frak g/\frak p\to \frak{su}(\Bbb V)/\frak q$.

\begin{lem*}
  In the situation of Theorem \ref{2.4}, the restriction of the
  curvature function $\underline{\ka}$ of $\underline{\om}$ to the
  domain of $\Ph$ is related to the curvature function $\ka$ of the
  conformal Cartan connection $\om$ by the commutative diagram 
$$
\xymatrix{%
\La^2(\frak g/\frak p) \ar[r]^{\ka\o \Ph}\ar[d] & \frak g\\
\La^2\big(\frak{su}(\Bbb V)/\frak q\big) \ar[r]^{\underline{\ka}} &
\frak{su}(\Bbb V),\ar[u]}
$$
where we view $\Cal G^{SU}$ as a subspace of $\Cal G$, and the
vertical arrows are the natural surjection, respectively  the
natural inclusion. This completely determines $\underline{\ka}$.
\end{lem*}
\begin{proof}
  Since $\om^{SU}$ is simply a restriction of $\om$, the definition of
  the curvature immediately implies that $K^{SU}$ is the restriction
  of $K$ to $\Cal G^{SU}$. (This restriction automatically has values
  in the algebra $\frak{su}(\Bbb V)\subset\frak g$.) In the proof of
  Theorem \ref{2.4} we have seen that on the domain of $\Ph$,
  $\underline{\om}$ coincides with $\Ph^*\om^{SU}$. Again by
  definition of the curvature, this shows that
  $\underline{K}=\Ph^*(K^{SU})$. Moreover, for $X\in\frak{su}(\Bbb V)$
  we get
  $T\Ph\o\underline{\om}^{-1}(X)=(\om^{SU})^{-1}(X)\o\Ph=\om^{-1}(X)\o\Ph$.
  Together with the above, we get for another element
  $Y\in\frak{su}(\Bbb V)$ and all points $a$ in the domain of $\Ph$
  the equation
  $\underline{K}(a)(\underline{\om}^{-1}(X),\underline{\om}^{-1}(Y))=
  K(\Ph(a))(\om^{-1}(X),\om^{-1}(Y))$. But the natural surjection
  $\frak g/\frak p\to \frak{su}(\Bbb V)/\frak q$ by definition maps
  $X+\frak p$ to $X+\frak q$ (and is even characterised by that), so
  we obtain the claimed relation between $\ka$ and $\underline{\ka}$.
  Since the domain of $\Ph$ contains $N\x\{e\}$ and the curvature
  function $\underline{\ka}$ is $Q$--equivariant, we see that this
  relation determines $\underline{\ka}$.
\end{proof}

\subsection{The characterisation theorem}\label{2.6}
To prove our main result, we have to verify one more property of the
conformal Killing field underlying a parallel complex structure on the
standard tractor bundle.

\begin{lem*}
  Let $(M,[g])$ be a conformal manifold endowed with a parallel
  complex structure $\Bbb J$ on the standard tractor bundle, let
  $\k^a$ be the conformal Killing field underlying $\Bbb J$ (see
  \ref{2.1}), and let $\Om_{ab}{}^A{}_B$ be the Cartan curvature. Then
  for the Levi--Civita connection $\nabla$ of a preferred scale in the
  sense of \cite[4.4]{fefferman}, we have
  $\Om_{ab}{}^A{}_B\nabla^a\k^b=0$.
\end{lem*}
\begin{proof}
  By the explicit formula for $\Om_{ab}{}^A{}_B$ in the proof of
  Proposition \ref{2.1} it suffices to show that
  $C_{abcd}\nabla^a\k^b=0$ and $A_{cab}\nabla^a\k^b=0$. In that proof,
  we have seen that $C_{abcd}\nabla^c\k^d=0$, so the first equation
  follows by the symmetries of the Weyl curvature. Using the covariant
  derivative of this equation, we obtain
  $$
  A_{cab}\nabla^a\k^b=\tfrac{1}{n-3}(\nabla^dC_{abdc})(\nabla^a\k^b)=
  \tfrac{-1}{n-3}C_{abdc}\nabla^d\nabla^a\k^b.
  $$
  There is an explicit formula for $\Bbb J_{AB}$ in part (5) of
  Proposition 4.4 of \cite{fefferman} in terms of any preferred scale.
  Expanding $0=Z^A{}_aZ^B{}_b\nabla_d\Bbb J_{AB}$ using this formula,
  we obtain
  $$
  0=\nabla_d\nabla_a\k_b+2\Rho_{d[a}\k_{b]}+2g_{d[a}\ell_{b]}.
  $$
  Inserting this into the above equation we get zero by
  the trace-freeness of $C_{abcd}$ and since $\k^aC_{abcd}=0$, which was
  observed in Proposition \ref{2.1}.
\end{proof}

Having this technical result at hand, we can proceed to the main
results for the characterisation.

\begin{prop*}
  The Cartan geometry $(N\x Q,\underline{\om})$ obtained in Theorem
  \ref{2.4} is torsion free (and hence regular) and normal. Thus it
  induces on the local leaf space $N$ a CR structure of hypersurface
  type, which is non--degenerate of signature $(p',q')$, as well as an
  appropriate root $\ce(1,0)$ of the canonical bundle, see
  \cite[2.3]{fefferman}.
\end{prop*}
\begin{proof}
  Lemma \ref{2.5} describes the curvature $\underline{\ka}$ of
  $\underline{\om}$. For $u\in\im(\Ph)\subset\Cal G^{SU}$, the map
  $\ka(u):\La^2(\frak g/\frak p ) \to \frak g$ descends to
  $\La^2(\frak{su}(\Bbb V)/\frak q)$ and has values in $\frak{su}(\Bbb
  V)$. First observe that since $\ka(u)$ is normal in the conformal
  sense, it is torsion free, so we conclude that the values actually
  lie in $\frak{su}(\Bbb V)\cap\frak p=\frak p^{SU}\subset\frak q$.
  This shows that $\underline{\om}$ is torsion free and hence regular
  as a Cartan connection on $N\x Q$. Regularity implies that the
  induced almost CR structure on $N$ (as described in Corollary
  \ref{2.4}) is non--degenerate of signature $(p',q')$ and partially
  integrable. Moreover, once we have proved that $\underline{\om}$ is
  normal, torsion freeness implies that the structure is integrable
  (and hence CR) by \cite[4.16]{Cap-Schichl}. The bundle $\ce(1,0)$
  can be defined as $(N\x Q)\x_Q L$ which immediately implies that it
  has the required properties.
  
  Hence the proof boils down to showing that, if we interpret $\ka(u)$
  as a map $\La^2(\frak{su}(\Bbb V)/\frak q)\to\frak{su}(\Bbb V)$, then
  it satisfies the normalisation condition $\partial^*\o\ka=0$ for
  Cartan geometries of type $(SU(\Bbb V),Q)$. Having proved this,
  normality of $\underline{\om}$ immediately follows from
  $Q$--equivariancy of $\partial^*$.
  
  To compute $\partial^*\o\ka$ we have to choose various bases.
  First, there is an abelian subalgebra $\frak p_+\subset\frak
  p\subset\frak g$ which consists of all maps $Z\in\frak g$ which
  vanish on $L_{\Bbb R}$ and map its real orthocomplement
  $L^{\perp_{\Bbb R}}_{\Bbb R}$ to $L_{\Bbb R}$.  Skew symmetry then
  implies that all values of $Z$ lie in $L^{\perp_{\Bbb R}}_{\Bbb R}$.
  The Killing form of $\frak g$ induces a duality between $\frak p_+$
  and $\frak g/\frak p$.
  
  Similarly, we obtain a subalgebra $\frak q_+\subset\frak
  q\subset\frak{su}(\Bbb V)$ by taking all maps $\underline{Z}$ which
  vanish on $L$ and map $L^\perp$ to $L$. Again, the Killing form
  induces a duality between $\frak{su}(\Bbb V)/\frak q$ and $\frak
  q_+$. However, in the complex case, there is a finer decomposition.
  Namely, there is a subspace $\frak q_2\subset\frak q_+$ of real
  dimension 1, which consists of those maps that vanish identically on
  $L^\perp$. (In the real case, the analogous condition already forces
  a map to vanish identically). The annihilator of $\frak q_2$ is the
  (real) codimension one subspace $\frak h\subset \frak{su}(\Bbb
  V)/\frak q$ used in the proof of Corollary \ref{2.4}.
  
  Since the Killing form on a simple Lie algebra is uniquely
  determined up to multiples by invariance, we may in both cases use
  the trace form on $\frak g$ for obtaining the dualities on both
  algebras. This only replaces $\partial^*$ by a non-zero multiple, so
  it has no effect on normality. Now let us chose a non-zero element
  $X_0\in\frak{su}(\Bbb V)$ such that $X_0+\frak q\notin\frak h$.
  Next, choose elements $X_1,\dots,X_{2n'}$ (where $n'=p'+q'$) such
  that $\{X_1+\frak q,\dots,X_{2n'}+\frak q\}$ is a basis of $\frak
  h$. Finally, put $X_{2n'+1}:=i\cdot\id\in\frak g$. Then $\{X_0+\frak
  p,\dots,X_{2n'+1}+\frak p\}$ is a basis for $\frak g/\frak p$. By
  duality, there are unique elements $Z_0,\dots,Z_{2n'+1}\in\frak p_+$
  such that $\tr(X_i\o Z_j)=\delta_{ij}$, so $\{Z_0,\dots,Z_{2n'+1}\}$
  is the dual basis of $\frak p_+$.
  
  On the other hand, for $j=0,\dots,2n'$ let $\underline{Z}_j$ be the
  complex linear part of the linear map $Z_j$, i.e.
  $\underline{Z}_j(v)=\tfrac{1}{2}(Z_j(v)-iZ_j(iv))$. This means that
  $Z_j-\underline{Z}_j$ is conjugate linear, and since each $X_i$ is
  complex linear, the map $X_i\o (Z_j-\underline{Z}_j)$ is conjugate
  linear.  But such a map has vanishing real trace, so we conclude
  that $\tr(X_i\o\underline{Z}_j)=\delta_{ij}$. Using this for
  $X_{2n'+1}=i\cdot\id$, we in particular see that the
  $\underline{Z}_j\in\frak{su}(\Bbb V)$ for $j=0,\dots,2n'$.
  
  Next, take a nonzero element $v\in L_{\Bbb R}$. By the definition of
  $\frak p_+$, each $Z_j(L^{\perp_{\Bbb R}}_{\Bbb R})\subset L_{\Bbb
  R}$, so $Z_j(iv)=av$ for some $a\in\Bbb R$.  But then $i\cdot\id\o
  Z_j$ maps $iv$ to $aiv$, and looking at an appropriate basis
  extending $\{v,iv\}$ one concludes that $\tr(i\cdot\id\o
  Z_j)=2a$. Hence $Z_j|_L=0$, and therefore $\underline{Z}_j|_L=0$ for
  $j=0,\cdots ,2n'$. For $w\in L^\perp$ all complex multiples of $w$
  lie in $L^{\perp_{\Bbb R}}_{\Bbb R}$, so $\underline{Z}_j(w)\in L$.
  Thus we conclude that $\underline{Z}_j\in\frak q_+$ for all
  $j=0,\dots,2n'$ and hence
  $\{\underline{Z}_0,\dots,\underline{Z}_{2n'}\}$ is the basis of
  $\frak q_+$ which is dual to the basis $\{X_0+\frak
  q,\dots,X_{2n'}+\frak q\}$ of $\frak{su}(\Bbb V)/\frak q$.
 
  Now the formula for $\partial^*(\ka(u)):\frak{su}(\Bbb V)/\frak
  q\to\frak{su}(\Bbb V)$ has two summands, see \cite[section
  5.1]{Yamaguchi}. The first summand is 
$$
X+\frak q\mapsto\sum_{j=0}^{2n'}[\ka(u)(X+\frak q,X_j+\frak
  q),\underline{Z}_j],
  $$
  where the bracket is the commutator of linear maps. Since
  $\ka(u)(X+\frak q,X_j+\frak q)$ is complex linear, this is the
  complex linear part of $\sum_{j=0}^{2n'}[\ka(u)(X+\frak q,X_j+\frak
  q),Z_j]$.  Now inside $\ka(u)$, we can replace $X+\frak q$ and
  $X_j+\frak q$ by $X+\frak p$ and $X_j+\frak p$, respectively, which
  just corresponds to viewing $\ka(u)$ as a map $\La^2\frak g/\frak
  p\to\frak{su}(\Bbb V)$.  Moreover, we may sum up to $2n'+1$ without
  changing the result since $\ka(u)$ descends to $\frak{su}(\Bbb
  V)/\frak q$. (This corresponds to the fact that $\k$ hooks trivially
  into $\ka$.) Thus have expressed the first summand as the complex
  linear part of $\sum_{j=0}^{2n'+1}[\ka(u)(X+\frak p,X_j+\frak
  p),Z_j]$, and the bracket is the same as in $\frak g$.  But by the
  choice of our bases, the condition that the conformal Cartan
  connection $\om$ is normal exactly reads as
  $\sum_{j=0}^{2n'+1}[\ka(u)(X+\frak p,X_j+\frak p),Z_j]=0$.
  
  Thus we are left with the second summand, which is
  \begin{equation*}\tag{$*$}
  X+\frak q\mapsto \tfrac{1}{2}
  \sum_{j=0}^{2n'}\ka(u)([X,\underline{Z}_j]+\frak q,X_j+\frak q).
  \end{equation*}
  We have already observed that $\underline{Z}_j$ vanishes on $L$, so
  $[X,\underline{Z}_j](L)=\underline{Z}_j(X(L))$. If $X+\frak
  q\in\frak h\subset\frak{su}(\Bbb V)/\frak q$, then
  $\underline{Z}_j(X(L))\subset\underline{Z}_j(L^\perp)\subset L$, and
  hence $[X,\underline{Z}_j]\in\frak q$. Therefore, it suffices to
  show that ($*$) vanishes for $X=X_0$.
  
  Take a nonzero element $v\in L_{\Bbb R}$. In the proof of Corollary
  \ref{2.4} we have seen that $X+\frak q\mapsto X(v)+L$ induces a
  linear isomorphism $\frak h\to L^\perp/L$, which was used to define
  the complex structure on $\frak h$. Fixing $w\in\Bbb V$ such that
  $\langle v,w\rangle=1$, one easily verifies, in a basis, that the
  restriction of the traceform to $\frak h\x\frak q_+$ is a nonzero
  real multiple of $(X+\frak q,\underline{Z})\mapsto\langle
  \underline{Z}(X(v)),w\rangle_{\Bbb R}$. 
  
  Now for our basis element $X_0$, by definition we have $X_0(v)\notin
  L^\perp$. Since $X_0(v)\neq 0$, it must be congruent to a nonzero,
  purely imaginary multiple of $w$ modulo $L^\perp$. Further, for
  $j=1,\dots,2n'$ we have $\underline{Z}_j\in\frak q_+$ so $\underline{Z}_j|_L=0$ and since 
  $\underline{Z}_j$ is skew Hermitian, this also implies $\underline{Z}_j(\Bbb V)\subset 
  L^\perp$. Hence $[X_0,\underline{Z}_j](v)=-\underline{Z}_j(X_0(v))\in L^\perp$. 
  On the one hand, this shows that $[X_0,\underline{Z}_j]+\frak q\in\frak
  h$. On the other hand, together with the above we see that there is
  a nonzero real number $a$ such that
  $$
  \langle X(v),i[X_0,\underline{Z}_j](v)\rangle_{\Bbb R}=-\langle
  X(v),i\underline{Z}_j(X_0(v))\rangle_{\Bbb R}=a\tr(X\o
  \underline{Z}_j).
  $$
  Now $(X+\frak q,Y+\frak q)\mapsto \langle X(v),Y(v)\rangle_{\Bbb
    R}$ defines a non--degenerate inner product on $\frak h$. Hence
  there are unique elements $Y_j+\frak q\in\frak h$ for
  $j=1,\dots,2n'$ such that $\langle X_k(v),Y_j(v)\rangle_{\Bbb
    R}=\delta_{kj}$. Applying the last displayed formula to $X=X_k$ we
  conclude that $[X_0,\underline{Z}_j]+q=-ia(Y_j+q)$. Therefore,
  normality is equivalent to
  \begin{equation*}\tag{$**$}
  \sum_{j=1}^{2n'}\ka(u)(i(Y_j+\frak q),X_j+\frak q)=0.
  \end{equation*}
  This expression admits an interpretation on $M$. Similar to the
  linear isomorphism $\frak h\to L^\perp/L$, a nonzero element
  $v\in\Bbb L_{\Bbb R}$ also gives rise to a linear isomorphism $\frak
  g/\frak p\to L^{\perp_{\Bbb R}}_{\Bbb R}/L_{\Bbb R}$. Now $\langle\ 
  ,\ \rangle_{\Bbb R}$ induces an inner product on $L^{\perp_{\Bbb
      R}}_{\Bbb R}/L_{\Bbb R}$ which can be carried over to $\frak
  g/\frak p$. Passing to the associated bundle $\Cal G\x_P(\frak
  g/\frak p)$ the resulting class of inner products induces the
  conformal structure on $M$. In particular, $L^\perp/L$ (together
  with the class of inner products induced by $\langle\ ,\ 
  \rangle_{\Bbb R}$) corresponds to a subquotient of the tangent
  spaces of $M$. Choosing a preferred metric from the conformal class,
  sections 4.4 and 4.5 of \cite{fefferman} show that this subquotient
  can be naturally identified with $\tilde
  H=\k^\perp\cap\ell^\perp\subset TM$. Moreover, by \cite[Proposition
  4.4]{fefferman}, the endomorphism $\nabla_a\k^b$ of $TM$ vanishes on
  $\k$ and $\ell$, and it represents the complex structure on $\tilde H$
  obtained from the identification with $L^\perp/L$. But this exactly
  means that the left hand side of ($**$) corresponds to
  $\Om_{abAB}\nabla^a\k^b$, which vanishes by the Lemma above.
\end{proof}

\begin{thm*}
  Let $(M,[g])$ be a conformal structure of signature $(2p'+1,2q'+1)$
  such that there is a parallel orthogonal complex structure $\Bbb J$
  on the standard tractor bundle. Then locally $M$ is conformally
  isometric to the Fefferman space of a CR manifold of hypersurface
  type, which is non--degenerate of signature $(p',q')$ endowed with
  an appropriate root $\ce(1,0)$ of the canonical bundle in such a way
  the $\Bbb J$ corresponds to the canonical complex structure
  introduced in \cite{fefferman}.
\end{thm*}
\begin{proof}
  Take an open subset $W\subset M$ and a local leaf space $\ps:W\to N$
  as in Theorem \ref{2.4}. By the Proposition above we get a CR
  structure of the required type on $N$, and $(N\x Q,\underline{\om})$
  is the canonical normal Cartan geometry associated to this
  structure. By construction, $P^{SU}\subset Q$ is the stabiliser of
  the real isotropic line $L_\Bbb R\subset L\subset\Bbb V$. Taking
  $(N\x Q)\x_Q L$ as $\ce(-1,0)$, the Fefferman space $\tilde N$ of
  $N$ is simply $N\x (Q/P^{SU})$, endowed with the Cartan geometry
  $(N\x Q\to N\x (Q/P^{SU}),\underline{\om})$. Theorem \ref{2.4}
  exactly says that this Cartan geometry is locally isomorphic to
  $(p:\Cal G^{SU}\to M,\om^{SU})$, which by the construction in
  \cite[Theorem 2.4 and 3.1]{fefferman} implies that $M$ is locally
  isometric to $\tilde N$ in a way compatible with the complex
  structures on standard tractors.
\end{proof}

\section{A strengthening of Sparling's characterisation}\label{3}

From some points of view, the local characterisation of Fefferman
spaces in Theorem \ref{2.6} is rather satisfactory. The
characterisation uses only conformally invariant data, and it is very
conceptual through its relation to conformal holonomy. On the other
hand, at first sight it would seem to be much weaker than the
 characterisation due to Sparling as
described in Theorem 3.1 of \cite{Graham}. This characterisation says
that if $(M,g)$ is a pseudo-Riemannian manifold, of appropriate
signature, which admits an isotropic Killing field $\k$ with the properties 
that it inserts
trivially into both the Weyl curvature and the Cotton tensor, while 
$Ric(\k,\k)>0$, then $M$ is locally
conformally isometric to a Fefferman space.

Now as shown in \cite[Proposition 2.2]{conf-Einst} and in
\cite[Corollary 3.5]{deformations}, parallel sections of the adjoint
tractor bundle $\Cal A$ of any conformal manifold $M$ are in bijective
correspondence with conformal Killing fields which insert trivially
into the Weyl curvature and the Cotton--York tensor. Given such a
conformal Killing field $\k$, one can view the corresponding parallel
adjoint tractor as an endomorphism $\Bbb J$ of the standard tractor
bundle $\Cal T$, and express the condition that $\Bbb J\o\Bbb J=-id$
in terms of $\k$ and its covariant derivatives.  Expressing $\Bbb
J\o\Bbb J=-id$ in this way, one is quickly led to see this includes
the properties required in Sparling's characterisation (see subsection
4.4 of \cite{fefferman}). However pushing through the calculation
na\"\i vely at first suggests that many more conditions arise than are
needed for the Sparling characterisation.  In fact a direct but
tedious calculation reveals that the remaining identities are all
differential consequences of the conditions in Sparling's
characterisation.

There is however an elegant and insightful way to do this using deeper
methods from conformal geometry. As a bonus we obtain a variant of
Sparling's characterisation which is stronger by dint of conformal
invariance.
\begin{thm*}
  Let $(M,[g])$ be a conformal manifold of dimension $n\geq 3$ and let
  $s\in\Ga(\Cal A)$ be a parallel section with underlying conformal Killing
  field $\k$, and suppose that $\k$ is isotropic. Then for any metric
  in the conformal class, with Levi--Civita connection $\nabla$ and
  Schouten tensor $\Rho$, the function
$$
\tfrac{1}{n^2}(\nabla_a\k^a)^2 -\k^a\Rho_{ab}\k^b -
\tfrac{1}{n}\k^a\nabla_a\nabla_b\k^b 
$$
is equal to some constant $\la\in\Bbb R$. This constant is
independent of the choice of metric from the conformal class 
and $s\o s=\la\id\in\Ga(\End(\Cal
T))$.
\end{thm*}
\begin{proof}
  As an endomorphism of $\Cal T$, $s$ is skew symmetric and hence $s\o
  s\in \Ga(S^2\Cal T)$. Now $S^2\Cal T$ invariantly decomposes into
  the trace-free part $S^2_0\Cal T$ and $\Bbb R\id$. Since $s$ is a
  parallel section of $\End(\Cal T)$, $s\o s$ is a parallel section of
  $S^2\Cal T$, and hence it is the sum of a parallel section of
  $S^2_0\Cal T$ and a constant multiple of the identity. Now the
  canonical filtration $\Cal T^1\subset\Cal T^0\subset\Cal T$ with
  $\Cal T^0:=(\Cal T^1)^\perp$ induces a filtration of $S^2_0\Cal T$.
  The largest filtration component is the image of $\Cal T^0\otimes
  \Cal T$ in $S^2_0\Cal T$ (under the obvious projection) and the
  quotient of $S^2_0\Cal T$ by this filtration component is isomorphic
  to $S^2(\Cal T/\Cal T^0)\cong\ce[2]$. If $t^{AB}$ is a section of
  $S^2_0\Cal T$ then the projection to this quotient is given by
  $X_At^{AB}X_B$.
  
  Now we can use the explicit formula for $s_{AB}$ in terms of
  $K_B:=Z_{Bb}k^b-\tfrac{1}{n}X_B\nabla_b\k^b$ from the proof of
  Proposition \ref{2.1}. This, in particular, implies $X^As_{AB}=K_B$,
  as well as 
$$
K^Bs_{BC}=-\tfrac{1}{n}K_C\nabla_b\k^b+\k^a\nabla_aK_C.
$$ This shows that $X^As_A{}^Cs_C{}^BX_B=-\k^a\k_a=0$. But a simple
corollary of the machinery of BGG sequences developed in
\cite{CSS-BGG} and \cite{David-Tammo} is that any parallel section of
an indecomposable tractor bundle can be recovered, by an invariant
differential operator, from its projection to the quotient by the
largest filtration component.

Hence $s\o s=\la\id_{\Cal T}$ for some constant $\la\in\Bbb R$, and we
can compute this constant for example as
$X^As_A{}^Cs_C{}^BY_B=K^Bs_{BC}Y^C$. The formula for $\la$ then
follows easily by expanding the above expression for $K^Bs_{BC}$.
\end{proof}

Using this, we obtain a strengthening of Sparling's characterisation as
well as surprising information about odd dimensional manifolds
admitting isotropic normal conformal Killing fields. 
\begin{cor*}
  Let $M$ be a pseudo-Riemannian manifold of dimension $n\geq 3$ 
and
  let $\k$ be an isotropic conformal Killing field, which inserts
  trivially into the Weyl curvature and (the 2--form indices of) the
  Cotton--York tensor. 

\noindent
(1) If $n$ is even and the constant
$$
\tfrac{1}{n^2}(\nabla_a\k^a)^2-\k^a\Rho_{ab}\k^b-
\tfrac{1}{n}\k^a\nabla_a\nabla_b\k^b 
$$
is negative, then $M$ is locally conformally isometric to a
Fefferman space $\tilde N$, of a CR manifold $N$, in such a way that
$\k$ is mapped to a generator of the vertical subbundle of $\tilde
N\to N$.

\noindent
(2) If $n$ is odd, then 
$$
\tfrac{1}{n^2}(\nabla_a\k^a)^2-\k^a\Rho_{ab}\k^b-
\tfrac{1}{n}\k^a\nabla_a\nabla_b\k^b=0
$$
and the adjoint tractor field $s$ corresponding to $\k$ satisfies
$s\o s=0$.
\end{cor*}
\begin{proof}
  (1) Passing to a constant rescaling of $\k$, we obtain a conformal
  Killing field such that the associated parallel adjoint tractor
  defines an almost complex structure. Then the result follows from
  the theorem.

\noindent
(2) Since $n$ is odd, the tractor bundle $\Cal T$ has odd rank, so the
skew symmetric map $s:\Cal T\to\Cal T$ has to be degenerate. Then $s\o
s=\la\id$ is only possible for $\la=0$, and the result follows.
\end{proof}

\subsection*{Remarks}
(1) If in part (1) of the corollary we assume that $\k$ is a Killing
field rather than just a conformal Killing field, then
$\nabla_a\k^a=0$, and the constant in question simply becomes
$-\k^a\Rho_{ab}\k^b$. Since $\k$ is isotropic, this is a nonzero
multiple of $Ric(\k,\k)$, and we recover Sparling's characterisation.
Note that in contrast to $Ric(\k,\k)>0$, negativity of the constant in (1)
does not evidently imply that $\k$ is nowhere vanishing. We get this as
a consequence of the corollary.

\noindent
(2) In section 4 of \cite{Graham} it is shown how to deduce a global
characterisation of Fefferman spaces from the local one.

\noindent
(3) In the odd--dimensional case, the kernel of $s$ is nontrivial in
each point. If the dimension of these kernels is constant, then they
form a smooth subbundle in the standard tractor bundle. Since $s$ is
parallel, this subbundle is invariant under the standard tractor
connection. In particular, this implies that the conformal holonomy
cannot act irreducibly.

\end{document}